\documentclass[12pt]{article}
\usepackage{amsmath,amsfonts,amssymb,amsthm}
\usepackage{fullpage}

\makeatletter
\def\section{\@startsection{section}{1}{\z@}{-3.5ex plus -1ex minus
  -.2ex}{2.3ex plus .2ex}{\large\bf}}
\def\subsection{\@startsection{subsection}{2}{\z@}{-3.25ex plus -1ex
  minus -.2ex}{1.5ex plus .2ex}{\normalsize\bf}}
\makeatother

\theoremstyle{definition}

\numberwithin{equation}{section}

\newcommand{\A}{\mathcal{A}}        
\renewcommand{\a}{\alpha}           
\newcommand{\B}{\mathcal{B}}        
\renewcommand{\b}{\beta}            
\newcommand{\braket}[2]{\langle#1\mathbin|#2\rangle} 
\newcommand{\C}{\mathbb{C}}         

\newcommand{\dn}{{\mathord{\downarrow}}} 
\DeclareMathOperator{\Dom}{Dom}     
\newcommand{\Dslash}{{D\mkern-11.5mu/\,}} 
\newcommand{\eps}{\varepsilon}      
\newcommand{\ev}{{\mathrm{ev}}}     
\newcommand{\ga}{\gamma}            
\renewcommand{\H}{\mathcal{H}}      
\newcommand{\half}{{\mathchoice{\oh}{\oh}{\shalf}{\shalf}}} 
\DeclareMathOperator{\ind}{Index}   
\newcommand{\K}{\mathcal{K}}        
\newcommand{\G}{\mathcal{G}}        
\DeclareMathOperator{\Mat}{Mat}     
\newcommand{\N}{\mathbb{N}}         
\newcommand{\nn}{\nonumber}         
\newcommand{\odd}{{\mathrm{odd}}}   
\newcommand{\oh}{{\tfrac{1}{2}}}    
\newcommand{\ooh}{{\tfrac{3}{2}}}   
\newcommand{\OP}{\mathrm{OP}}       
\newcommand{\otto}{\leftrightarrow} 
\newcommand{\ox}{\otimes}           
\newcommand{\piappr}{\underline{\pi}} 
\newcommand{\R}{\mathbb{R}}         
\DeclareMathOperator{\Res}{Res}     
\newcommand{\sepword}[1]{\quad\mbox{#1}\quad} 
\newcommand{\Sf}{\mathbb{S}}        
\newcommand{\sg}{\sigma}            
\newcommand{\shalf}{{\scriptstyle\frac{1}{2}}} 
\newcommand{\ssesq}{{\scriptstyle\frac{3}{2}}} 
\newcommand{\sesq}{{\mathchoice{\ooh}{\ooh}{\ssesq}{\ssesq}}} 
\DeclareMathOperator{\Tr}{Tr}       
\newcommand{\ul}[1]{\underline{#1}} 
\newcommand{\up}{{\mathord{\uparrow}}} 
\newcommand{\Z}{\mathbb{Z}}         
\renewcommand{\:}{\colon}           
\def\<#1,#2>{\langle#1,#2\rangle}   
\def\CD{{\mathcal D}}
\def\DO{{\mathcal D}}

\newbox\ncintdbox \newbox\ncinttbox
\setbox0=\hbox{$-$}
\setbox2=\hbox{$\displaystyle\int$}
\setbox\ncintdbox=\hbox{\rlap{\hbox
    to \wd2{\kern-.1em\box2\relax\hfil}}\box0\kern.1em}
\setbox0=\hbox{$\vcenter{\hrule width 4pt}$}
\setbox2=\hbox{$\textstyle\int$}
\setbox\ncinttbox=\hbox{\rlap{\hbox
    to \wd2{\kern-.05em\box2\relax\hfil}}\box0\kern.1em}
\newcommand{\ncint}{\mathop{\mathchoice{\copy\ncintdbox}%
    {\copy\ncinttbox}{\copy\ncinttbox}%
    {\copy\ncinttbox}}\nolimits}

\begin{document}

\title{{\bf
The local index formula for quantum \boldmath $SU(2)$}
}
\author{~\\
{\bf
Ludwik D\c{a}browski}\\
\normalsize
~\\
Scuola Internazionale Superiore di Studi Avanzati,\\
\normalsize
Via Beirut 2-4, 34014 Trieste, Italy
}
\date{}

\maketitle

\thispagestyle{empty}

~\\
\begin{abstract}
\noindent
The local index formula of Connes--Moscovici for the
isospectral noncommutative geometry recently constructed
on quantum $SU(2)$ \cite{Naiad,DLSSV2} is discussed.
The cosphere bundle and the dimension spectrum as well as the local 
cyclic cocycles yielding the index formula, are presented.\\
~\\
Talk at the Workshop "Traces in Geometry, 
Number Theory and Quantum Fields", 
MPIM Bonn, October 24-28, 2005.
\end{abstract}

\vspace{1pc}



\section{Introduction}

In noncommutative geometry a` la Connes \cite{Book} 
{\em Riemannian spin geometry} is described by {\em  spectral triple} $(\A,\H,D)$,
which consists of a (suitable) $*$-algebra $\A$, 
represented by bounded operators on a Hilbert space $\H$ 
and a self-adjoint operator $D$ on $\H$ with compact resolvent 
and bounded commutators with $\A$.
The canonical spectral triple 
associated to a spin manifold $M$ with a Riemannian metric $g$ is
$(C^\infty (M, \C), L^2(\sigma , vol_g), \Dslash )$, 
where $\sigma $ is 
the Dirac spinor bundle over $M$ and $\Dslash $ is the Dirac operator 
of the Levi-Civita connection of $g$.
This classical spectral triple satisfies seven additional conditions 
which has been formulated and postulated also for noncommutative algebras
$\A$ \cite{ConnesGrav}. One of them, \textit{regularity} (or \textit{smoothness}),
 permits to introduce \cite{ConnesMIndex}
the pseudodifferential calculus.
Another one, {\it dimension}, allows to define generalized zeta functions
and dimension spectrum, which are tools for the local index theorem 
of Connes-Moscovici \cite{ConnesMIndex},
 a powerful algorithm for performing complicated local computations 
by neglecting plethora of irrelevant details.
This occurs because these formulae employ exotic traces given in terms of residues
rather than the usual traces (which e.g. in position representation 
require multiple integrals over the whole $M$) employed in the index theorem
formulated in terms of Fredholm modules.
\\
~\\
%
An important area to implement and probe these ideas is quantum groups,
to start with the best known $SU_q(2)$ (the quantum $SU(2)$).
On $SU_q(2)$ a preliminary candidate \cite{BibikovK} for the operator $D$ 
had unbounded commutators with  $\A$ \cite{Goswami}.
The first spectral triple, constructed in~\cite{ChakrabortyPEqvt},
was `singular' (in the sense that it has no $ \lim_{q\to 1}$).
It has been extensively studied in \cite{ConnesSUq}
using the concept of a (quantum) cosphere bundle
$\Sf_q^*$ on $SU_q(2)$, that considerably simplifies the
computations of local index formula by 
removing the irrelevant smoothing operators 
(which give no contribution to the residues appearing in the
local cyclic cocycle).\\
~\\
I will present here an interesting isospectral 
bi-equivariant $3^+$-summable spectral triple 
recently constructed in \cite{Naiad} and analyzed  in \cite{DLSSV2} 
along the lines of \cite{ConnesSUq}.
The resulting cosphere bundle coincides with that in \cite{ConnesSUq},
as well as the dimension spectrum $\Sigma$, given by the set $\{1,2,3\}$.
As computed in \cite{MasudaNW} the cyclic cohomology of the algebra 
$\A(SU_q(2))$ is given in terms of a single generator. 
This element can be expressed in terms of a single local cocycle,
whose form in \cite{DLSSV2} is slightly different from the one in
\cite{ConnesSUq}.
A nice exemplification of the theory is computation of the index of $D$ 
coupled with the unitary representative of the generator of $K_1(\A)$.
The associated Fredholm module turns out to be 1-summable.

\subsection{Preliminaries on the pseudodifferential calculus}

The \textit{regularity} (or \textit{smoothness}) requirement for a spectral 
triple is that $\A\cup [D,\A] \subset \bigcap_{n=1}^\infty \Dom \delta^n$,
where 
$\delta(T) := |D|T - T|D|$.
This condition permits to introduce \cite{ConnesMIndex}
the pseudodifferential calculus as follows.\\
First, let $\H^s := \Dom(D^s)$ for $s \in \R$
 be the analogue of Sobolev spaces and let 
$\H^\infty := \bigcap_{s\geq 0} H^s$.
(Assume for simplicity that that $|D|$ is invertible; 
for the noninvertible see e.g. \cite{CareyPRS} 
and that $\H^\infty$ is a core for~$|D|$).\\
An operator $T\: \H^\infty \to \H^\infty$ such that 
$
|D|^{-\a} T \in \bigcap_{n=1}^\infty \Dom \delta^n \ ,
$
for $\a\in \R$, is said to have order $\a$.
(Such $T$ automatically extends to a bounded operator 
from $\H^{\a +s}$ to $\H^s$ for all $s \geq 0$). 
Let  $\OP^\a$ denote the set of operators of order $\leq\a$.
In particular, $\OP^0 = \bigcap_{n=1}^\infty \Dom \delta^n$, 
the algebra of operators of order~$\leq 0$ includes 
$\A \cup [D,\A]$ and their iterated commutators with~$|D|$. 
Moreover, $[D^2, \OP^\a] \subset \OP^{\a+1}$
and $\OP^{-\infty} := \bigcap_{\a\leq 0} \OP^\a$ 
is a two-sided ideal in~$\OP^0$.\\
With this set up, the algebra $\CD = \cup \CD_k$
of {\it differential operators} is just the smallest algebra of operators 
on $\H^\infty$ generated by $\A\cup [D, \A]$ and 
filtered by the order $k\in\Bbb N$ in such a way that
$[D^2, \DO_k] \subset \DO_{k+1}$.
Then the algebra $\Psi$ of {\it pseudodifferential operators},
is generated by operators, which modulo $OP^{\a}$ for any $\a\in\R$ 
are of the form $ T\, D^{-n}$ for some $n$ and some $T\in \CD$.
In particular $\Psi^0$ of order $\leq 0$ is the algebra generated by 
$\bigcup_k \delta^k( \A\cup [D, \A])$.\\
The algebra structure on $\Psi$ can be read off in terms of an
\textit{asymptotic expansion}: $T \sim \sum_{j=0}^\infty T_j$
whenever $T$ and each $T_j$ are operators from $\H^\infty$ to
$\H^\infty$; and for each $m \in \Z$, there exists $N$ such that for
all $M > N$, the operator $T - \sum_{j=1}^M T_j$ has order
$\leq m$. For instance, for complex powers of $|D|$ 
(e.g. defined by the Cauchy formula) there is a binomial expansion:
$$
[|D|^z, T] \sim \sum_{k=1}^\infty 
(k\, !)^{-1} z(z\! -\! 1)\dots (z\! -\! k\! +\! 1)\,
\delta^k(T)\,|D|^{z-k} .
$$
Another postulated requirement is that of {\it dimension}:
$\exists$ $n\in \N$ s.t. 
the eigenvalues (with multiplicity) of $|D|^{-n}$, 
$\mu_k = O(k^{-1})$ as $k\to \infty$.
Then
for $k > n$,
$D^{-k}$ is 
trace-class and 
the ~{\it dimension spectrum} $\Sigma$  is the set 
of the singularities of zeta functions
\begin{equation}
\zeta_\b (z) = {\rm Trace}_{\cal H} (\b\, |D|^{-z}),
~~\forall \b\in\Psi^0\ .
\end{equation}
Assuming $\Sigma$ to be discrete with simple poles only 
the Wodzicki-type residue functional
$$
\int  T |D|^{-n} := Res_{z=n} {\rm Trace}(T|D|^{-z})$$
is tracial on $T\in \Psi$.
%

\section{The isospectral geometry of \boldmath $SU_q(2)$}
\label{sec:iso-ge}

\subsection{Spectral triple}
The spectral triple $(\A,\H,D)$ of \cite{Naiad} 
can be written in the following form:\\
~\\
\hspace*{-.30cm}$\bullet$ 
$\A $ is the $*$-algebra \cite{W} generated by $a$ and~$b$ with 
$$
ba = q ab,  \quad  b^*a = qab^*, \quad bb^* = b^*b,\quad
a^*a + q^2 b^*b = 1,  \quad  aa^* + bb^* = 1.
$$
(Here $0 < q < 1$ and $a \otto a^*$, $b \otto -b$ are exchanged
with respect to \cite{ChakrabortyPEqvt} and~\cite{ConnesSUq}).\\
~\\
\hspace*{-.30cm}$\bullet$ 
The Hilbert space of spinors $\H $
has an orthonormal basis $v_{x,y,s}^j$ where
\begin{equation}
j = 0, \half, 1, \dots; ~~
x = 0, 1, \dots, 2j; ~~y = 0, 1, \dots, 2j\! +\! 1; ~~s = \up, \dn ;
\label{eq:new-basis}
\end{equation}
with the convention that the $\dn$ component is zero if $y = 2j$ or $2j+1$.\\
%
The spinor representation $\pi$ of $\A$ reads
$\pi(a) := a_+ + a_-$,~ $\pi(b) := b_+ + b_-$, where 
\begin{align}
a_+ v_{x,y,\up}^j 
&= q^{(x+y-2j-1)/2} [x\! +\! 1]^\half 
\left(
    \frac{q^{-j -1/2}[y\! +\! 1]^{1/2}}{[2j+2]} v_{x+1,y+1,\up}^{j+\half} +
    \frac{q^{\half}[2j\! -\! y\! +\! 1]^{\half}}{[2j+1]\,[2j+2]}  v_{x+1,y,\dn}^{j+\half}
\right),\nonumber\\
a_+  v_{x,y,\dn}^j
&= q^{(x+y-2j-1)/2} [x\! +\! 1]^\half 
    \frac{q^{-j}[y\! +\! 1]^{\half}}{[2j+1]} v_{x+1,y+1,\dn}^{j+\half}\ ,\nonumber\\
a_-  v_{x,y,\up}^j
&= q^{(x+y-2j-1)/2} [2j\! -\! x]^\half 
    \frac{q^{j+1}[2j\! -\! y\! + \! 1]^{\half}}{[2j+1]} v_{x,y,\up}^{j-\half}\ ,
\nonumber\\
a_- v_{x,y,\dn}^j 
&= q^{(x+y-2j-1)/2} [2j\! -\! x]^\half
\left(
    -\frac{q^{\half}[y\! +\! 1]^{\half}}{[2j][2j+1]} v_{x,y+1,\up}^{j-\half} +
    \frac{q^{j+\half}[2j\! -\! y\! -\! 1]^{\half}}{[2j]}  v_{x,y,\dn}^{j-\half}
\right),\nonumber\\
b_+ v_{x,y,\up}^j 
&= q^{(x+y-2j-1)/2} [x\! +\! 1]^\half 
\left(
    \frac{[2j\!- \! y \! +\! 2]^{\half}}{[2j+2]} v_{x+1,y,\up}^{j+\half} -
    \frac{q^{-j-1}[y]^{\half}}{[2j+1]\,[2j+2]}  v_{x+1,y-1,\dn}^{j+\half}
\right),
\nonumber\\
b_+  v_{x,y,\dn}^j
&= q^{(x+y-2j-1)/2} [x\! +\! 1]^\half 
 \cdot
    \frac{q^{-\half}[2j\! -\! y]^{\half}}{[2j+1]} v_{x+1,y,\dn}^{j+\half}\ ,
\nonumber\\
b_- v_{x,y,\up}^j
&= -q^{(x+y-2j-1)/2} [2j\! -\! x]^\half 
 \cdot
    \frac{q^{-\half}[y]^{\half}}{[2j+1]} v_{x,y-1,\up}^{j-\half}\ ,
\nonumber\\
b_- v_{x,y,\dn}^j 
&= -q^{(x+y-2j-1)/2} [2j\! -\! x]^\half
\left(
    -\frac{q^{j}[2j\! -\! y]^{\half}}{[2j][2j+1]} v_{x,y,\up}^{j-\half} +
    \frac{[y]^{\half}}{[2j]}  v_{x,y-1,\dn}^{j-\half}
\right).
\nonumber
\label{eq:spin-rep}\nonumber
\end{align}
\noindent
(Here $[N] := (q^{-N} - q^N)/(q^{-1} - q)$).\\
This representation is $U_q(su(2))\otimes U_q(su(2))$-equivariant
and unique (up to phases) on $\H$.\\
~\\
\hspace*{-.3cm}$\bullet$ 
The Dirac operator 
\begin{equation}
D v_{x,y,\up}^j = (2j\! +\! \sesq )\, v_{x,y,\up}^j, ~~
D v_{x,y,\dn}^j = - (2j\! +\! \half )\,  v_{x,y,\dn}^j, ~~
\label{eq:classical-evs}
\end{equation}
whose spectrum (with multiplicity) coincides with that of the
classical Dirac operator on the round sphere $S^3$,
is one of a family of operators in \cite{Naiad} 
($U_q(su(2))\otimes U_q(su(2))$-invariant, assymptotically diagonal
 with linear spectrum satisfying a modification of {\em reality}
 and {\em first order condition}).

\subsection{Analytic properties}
We shall use the polar decomposition $D = F\,|D|$, where
$|D|:= (D^2)^\half$ and $F = D/|D|$, and the orthogonal projectors
 $P^\up := \half(1 + F)$, $P^\dn := \half(1 - F) = 1 - P^\up$,
whose range spaces are respectively $\H^\up$, spanned by $v_{x,y,\up}^j$
  and $\H^\dn$ spanned by $v_{x,y,\dn}^j$.
Let $\delta^n (T)= [~|D|, \dots [~|D|, T\underbrace{]\dots ]}_n$.\\
The explicit computation in \cite{DLSSV2} shows that
\begin{equation} \delta(a_+) =
P^\up a_+ P^\up + P^\dn a_+ P^\dn, \qquad \delta(a_-) = - P^\up a_-
P^\up - P^\dn a_- P^\dn.  
\label{eq:delta-apm}  \nonumber
\end{equation}
Hence
$\delta(\pi(a)) =
\delta(a_+) + \delta(a_-)$ is bounded.
Next, using 
$\delta([D,\cdot ]) = [D,\delta(\cdot )]$,
\begin{equation} \delta([D,a_+]) = P^\up a_+ P^\up - P^\dn a_+ P^\dn,
\qquad \delta([D,a_-]) = P^\up a_- P^\up - P^\dn a_- P^\dn.
\label{eq:delta-Dapm} \nonumber
\end{equation} 
Calculations for $b$ give similar results and thus 
$\A\cup [D,\A]$ is contained in $ \Dom \delta$,
and by iteration, also in $ \Dom \delta^k$ for $k\in\N$.  
Therefore, 
the triple $(\A,\H,D)$ is a 
$U_q(su(2))\otimes U_q(su(2))$-equivariant
regular $3^+$-summable spectral triple.
(This holds also for a suitable completion of $\A(SU_q(2))$).\\ 
~\\
The details of the calculations above show that 
$\Psi^0$, the algebra spanned by $\delta^k(\A)$ and $\delta^k([D,\A])$,
for all $k \geq 0$, is in fact generated by the diagonal-corner operators
$P^\up a_\pm P^\up$, $P^\dn a_\pm P^\dn$, $P^\up b_\pm P^\up$, $P^\dn
b_\pm P^\dn$ together with the other-corner operators 
$P^\dn a_+ P^\up$, $P^\up a_- P^\dn$, $P^\dn b_+ P^\up$, $P^\up b_- P^\dn$.\\
The algebra $\B$ spanned by all $\delta^k(\A)$ for $k \geq 0$ 
is generated by the diagonal operators 
\begin{equation} 
\tilde a_\pm := \pm \delta(a_\pm) = P^\up a_\pm P^\up + P^\dn a_\pm P^\dn,
\label{gensB}
\end{equation} 
\begin{equation} 
\tilde b_\pm := \pm \delta(b_\pm) = P^\up b_\pm P^\up + P^\dn
b_\pm P^\dn, \label{eq:extended-pi} 
\nonumber
\end{equation} 
and by (off-diagonal) operators 
\begin{equation} 
\tilde a_{/} = P^\dn a_+ P^\up + P^\up a_- P^\dn, ~
\tilde b_{/} = P^\dn b_+ P^\up + P^\up b_- P^\dn\ .
\label{gensB/}
\end{equation}
%
Note that $\Psi^0$ is generated by (its subalgebra) $\B$ and by $P^\up$.\\
~\\
\noindent
In the sequel we shall need operators 
$\piappr(a) := \ul{a}_+ + \ul{a}_-$ and
$\piappr(b) := \ul{b}_+ + \ul{b}_-$, where 
\begin{align}
\ul{a}_+ v_{x,y,s}^j
&= \sqrt{1 - q^{2x+2}} \sqrt{1 - q^{2y+2}} \, v_{x+1,y+1,s}^{j^+}\ ,
\nn \\
\ul{a}_- v_{x,y,s}^j &= q^{x+y+1} \, v_{x,y,s}^{j^-}\ ,
\nn \\
\ul{b}_+ v_{x,y,s}^j &= q^y \sqrt{1 - q^{2x+2}} \, v_{x+1,y,s}^{j^+}\ ,
\nn \\
\ul{b}_- v_{x,y,s}^j &= - q^x \sqrt{1 - q^{2y}} \, v_{x,y-1,s}^{j^-}\ .
\label{eq:approx-repn-bis}
\end{align}
These formulae coincide with those in \cite[Sec.~6]{ConnesSUq}
up to the exchange $a \otto a^*$, $b \otto -b$ and 
a doubling of the Hilbert space ($s=\up,\dn$).
Using a truncation coming from 
\begin{align*}
\bigl((q^{-1} - q)[n] \bigr)^{-1} - q^n &= q^{3n} + O(q^{5n}),
\\
1 - \sqrt{1 - q^\a} & \leq q^\a,  \qquad{\rm for ~any} \ \a \geq 0;
\end{align*}
it can be seen that, for $x=a$ or $x=b$, ~$\piappr(x)$ approximate $\pi(x)$
up to operators (given by matrices) of rapid decay and so belonging 
to $\OP^{-\infty}$.
Hence, $\piappr(x)$ can be used instead of $\pi(x)$
when dealing with the local cocycle 
in the local index theorem in the sequel.
\\
Moreover, the operators $\ul{x}_{\pm}$, for $x= a,b$,
satisfy simple commutation rules
$$
[|D|, ~\ul{x}_{\pm}] = \pm\ul{x}_{\pm}, \quad
[D, ~\ul{x}_{\pm}] = \pm F\, \ul{x}_{\pm}, \quad
[F, \ul{x}_\pm ]= 0\ .
$$

\subsection{The cosphere bundle}
\label{sec:cosphere}

The `cosphere bundle' of $SU_q(2)$ constructed in \cite{ConnesSUq} using the
regular representation of $\A$ and the one obtained in \cite{DLSSV2} from 
the spinor representation are isomorphic.
To see this let $\pi_\pm$ be two known (bounded) representations
of $\A(SU_q(2))$ on the Hilbert space
$\ell^2(\N)$ with the standard orthonormal basis $\eps_k$, $ k \in \N$,
determined by
\begin{equation}
\pi_\pm(a) \,\eps_k := \sqrt{1-q^{2k+2}} \,\eps_{k+1},  \qquad
\pi_\pm(b) \,\eps_k := \pm q^k \,\eps_k.
\label{eq:pi-pm}
\end{equation}
By sending 
\begin{equation}
v_{x,y,s}^j \mapsto \eps_{j,x,y,s}:=
\eps_x \otimes \eps_y \otimes \eps_j\otimes \eps_s
\label{eq:identific}
\end{equation}
we identify the Hilbert space $\H$ 
with the subspace $\H'\subset
\ell^2(\N)_x \ox \ell^2(\N)_y \ox \ell^2(\Z)_{2j} \ox \C^2$
given by the restrictions of indices \eqref{eq:new-basis}.\\
This yields the correspondence
\begin{align}
\ul{a}_+ &\otto \pi_+(a) \ox \pi_-(a) \ox V \ox 1_2,
\nn \\
\ul{a}_- &\otto - q\,\pi_+(b) \ox \pi_-(b^*) \ox V^* \ox 1_2,
\nn \\
\ul{b}_+ &\otto - \pi_+(a) \ox \pi_-(b) \ox V \ox 1_2,
\label{eq:corr-smoothing}\\
\ul{b}_- &\otto - \pi_+(b) \ox \pi_-(a^*) \ox V^* \ox 1_2,
\nn
\end{align}
where $V$ is the unilateral shift operator
$\eps_{2j} \mapsto \eps_{2j+1}$ in $\ell^2(\Z)$.\\
~\\
Few remarks on \eqref{eq:corr-smoothing} are in order:\\
- Note that the two first factors (taking into consideraton the overall 
minus signs in the lats two equations) reproduce  
the Hopf tensor product of the representations $\pi_-$ and $\pi_+$,
with respect to the standard coproduct of $SU_q(2)$.\\
- The shift $V$ is best 
encoded using the $\Z$-grading due to the one-parameter group $\ga$
of automorphisms 
(playing the role of `geodesic flow', see \cite{ConnesSpec}), where
$\ga(t) :  T \mapsto  e^{it|D|} T e^{-it|D|}$
for any operator $T$ on $\H$.
On the subalgebra of ``diagonal''
operators $T = P^\up T P^\up + P^\dn T P^\dn$, 
$\ga$ detects the correct shift of $j$, for example, 
$\ga(t) :\ul{x}_{\pm}  \mapsto e^{\pm it}\ul{x}_{\pm}$, when $x=a, b$.\\
- \eqref{eq:corr-smoothing} coincides with \cite[(204)]{ConnesSUq}
apart from the last factor $1_2$ (and different conventions).\\
~\\
In the rest of the section the meaning of the correspondence 
\eqref{eq:corr-smoothing} will be clarified.\\
Since $b - b^* \in \ker \pi_\pm$, the reps $\pi_\pm$ are not faithful.
Let $\A(D^2_{q\pm})$ be the two quotient algebras defined by
\begin{equation}
0 \to \ker \pi_\pm \to \A(SU_q(2))
\xrightarrow{r_\pm} \A(D^2_{q\pm}) \to 0.
\label{eq:q-disks}
\end{equation}
In $\A(D_{q\pm}^2)$ one has
(omitting the quotient maps $r_\pm$)
\begin{gather}
b = b^*, \qquad b a = q\, a b,  \qquad    a^*b = q\, b a^*,
\nn \\
a^*a + q^2 b^2 = 1,  \qquad  aa^* + b^2 = 1.
\label{eq:Sq2-relns}
\end{gather}
These are just relations of the equatorial Podle\'s sphere 
$S^2_q$~\cite{Podles} (modulo $b \mapsto q^{-1} b$ plus $q \mapsto q^2$). 
But the spectrum of $\pi_\pm(b)$ being $\pm$tive,
$\A(D^2_{q\pm})$ actually describe the two hemispheres of $S^2_q$ 
(thought of as quantum disks).\\
~\\
In \cite{DLSSV2} it was proven that there exists a $*$-homomorphism
\begin{equation}
\rho: \B \to \A(D^2_{q+}) \ox \A(D^2_{q-}) \ox \A(S^1)
\label{eq:symbol-map-bis}
\end{equation}
defined on the generators (\ref{gensB},\ref{gensB/}) by
\begin{align*}
\rho(\tilde a_+) &:= r_+(a) \ox r_-(a) \ox u,  &
\rho(\tilde a_-) &:= -q\,r_+(b) \ox r_-(b^*) \ox u^* ,
\\
\rho(\tilde b_+) &:= - r_+(a) \ox r_-(b) \ox u, &
\rho(\tilde b_-) &:= - r_+(b) \ox r_-(a^*) \ox u^* ,
\end{align*}
\begin{equation}
\rho(\tilde a_{/}) = 0 = \rho(\tilde b_{/})\ .
\end{equation}
%
To see this, 
since the factors $u, u^*$ take care of the $j$-dependence,
it suffices to show that
\begin{align*}
\rho_\bullet (\tilde a_+) &:= \pi_+(a) \ox \pi_-(a), &
\rho_\bullet (\tilde a_-) &:= -q\,\pi_+(b) \ox \pi_-(b^*),
\\
\rho_\bullet (\tilde b_+) &:= - \pi_+(a) \ox \pi_-(b), &
\rho_\bullet (\tilde b_-) &:= - \pi_+(b) \ox \pi_-(a^*),
\end{align*}
\begin{equation}
\label{ker}
\rho_\bullet (\tilde a_{/}) = 0 = \rho_\bullet (\tilde b_{/}) \ ,
\end{equation}
determine a 
$*$-homo\-morph\-ism ~ $\rho_\bullet : \B \to \A(D^2_{q+}) \ox \A(D^2_{q-})$.\\
(Here $r_\pm(x)$ for $x\in \A(SU_q(2))$, is identified with its faithful
representant $\pi_\pm(x)$.)
~\\
Let $\Pi: \H \to \ell^2(\N) \otimes \ell^2(\N)$
be (unbounded) operator 
$ v_{x,y,s}^j \mapsto \eps_{x,y}:= \eps_x \otimes \eps_y$
(it `forgets' the index $j$ and $s$).\\
Define the map $\rho_\bullet $ by associating  to $T \in \B$ 
the operator $\rho_\bullet (T)$ on $\ell^2(\N) \otimes \ell^2(\N)$,
given by
\begin{equation}
\rho_\bullet (T)
\, \eps_{x,y} = \lim_{j \to \infty} \Pi(T v^j_{x,y,s}).
\end{equation}
It is well-defined, since $T$ is polynomial in generators of $\B$ 
that are weighted shifts (uniformly bounded) in $x,y,j$.\\
Now, using estimates of \cite{Naiad}, it can be directly verified
that $\rho_\bullet$ vanishes on the (rapidly decreasing) operators 
$\tilde a_{/}$, $\tilde b_{/}$, $\ul a_\pm - \tilde a_\pm$
and $\tilde b_\pm - \ul b_\pm$.\\
Hence the last eq. of \eqref{ker} holds.
Moreover $\tilde a_\pm$ and $\tilde b_\pm$ 
can be replaced by respectively $\ul a_\pm$ and $\ul b_\pm$ 
in the first four equations,
which are then satisfied because the coefficients in 
$\ul a_\pm$ and $\ul b_\pm$ (c.f. \eqref{eq:approx-repn-bis}) 
are $j$-independent, e.g.
$$
\rho_\bullet (\tilde a_+) \eps_{x,y}  =\rho_\bullet (\ul a_+) \eps_{x,y} =
\lim_{j \to \infty}  \sqrt{1-q^{2x+2}} \sqrt{1-q^{2y+2}} \Pi(v_{x+1, y+1,s}^{j^+})\\
$$
$$
=  \sqrt{1-q^{2x+2}} \sqrt{1-q^{2y+2}} \eps_{x+1, y+1} = 
(\pi_+ (a) \otimes \pi_-(a) ) \eps_{x,y}.
$$
Moreover, since the product of $\ul a_\pm$ and $\ul b_\pm$ still does not contain 
$j$-dependent coefficients, $\rho_\bullet $ respects the multiplication in $\B$.
Finally, $\rho_\bullet $ is an algebra homomorphism
by linearity of  $\lim$.\\
~\\
The range of the map $\rho$ in $\A(D^2_{q+}) \ox A(D^2_{q-}) \ox \A(S^1)$,
denoted $\A(\Sf_q^*)$, is called (algebra of) {\it cosphere bundle} on $SU_q(2)$
(and $\rho$ is `symbol map').\\
It deserves its name since it corresponds to complete symbols,
i.e. scalar pseudodifferential operators modulo the smoothing operators,
which in our case coincide with diagonal pseudodifferential operators 
modulo the smoothing operators\\
It should be mentioned that $\Sf_q^*$ for our geometry coincides with 
the cosphere bundle constructed in \cite{ConnesSUq} and 
the symbol map $\rho$ rectifies the correspondence \eqref{eq:corr-smoothing}.
It can be extended to $\Psi^0$ by setting $\rho (P) = P$
and to appropriate smooth algebras.\\  
~\\
Denote by $Q$ the orthogonal projector on
$\ell^2(\N)\ox \ell^2(\N) \ox \ell^2(\Z) \ox \C^2$ with range $\H'$
(the Hilbert subspace identified with~$\H$).
Then \eqref{eq:corr-smoothing} in
combination with \eqref{eq:symbol-map-bis}
implies that for all $T \in \B$
\begin{equation}
T - Q (\rho(T) \ox 1_2) Q \in \OP^{-\infty}\ ,
\label{eq:symbol-smoothing}
\end{equation}
where the first $T$ is viewed as an operator on $\H'$ via \eqref{eq:identific}
and $\rho(T)$ acts on 
$\ell^2(\Z)$ via Fourier transform from ~$S^1$ to $\Z$.

\subsection{The dimension spectrum}
Let $\tau_1$ and
$\tau_0^\up$, $\tau_0^\dn$ be three functionals 
on $\A(D_{q\pm}^2)$ defined by
\begin{align*}
\tau_1(x) &:= \frac{1}{2\pi} \int_{S^1} \sg(x),
\\
\tau_0^\up(x) &:= \lim_{N\to\infty} \left(
\Tr_N \pi_\pm(x) - (N+\sesq) \tau_1(x)  \right),
\\
\tau_0^\dn(x) &:= \lim_{N\to\infty} \left(
\Tr_N \pi_\pm(x) - (N+\half) \tau_1(x)\right),
\end{align*}
where\\
- $\sg \: \A(D^2_{q\pm}) \to \A(S^1)$ is the $*$-homomorphism 
(known as `symbol map')
that includes 
$S^1 = \partial D^2_{q\pm} $
(the equator of $S^2_q$) into $D^2_{q\pm} $
\begin{align}
\sg(r_\pm(a)) := u; \qquad \sg(r_\pm(b)) := 0,
\label{eq:symbol-map}
\end{align}
with $u$ being the unitary generator of $\A(S^1)$;\\
- $\Tr_N$ is the truncated trace
$$
\Tr_N(T) := \sum_{k=0}^N \braket{\eps_k}{T\eps_k}
$$
- the constants $\sesq$ and $\half$
are chosen to simplify the residues in the sequel.\\
~\\
A straightforward calculation shows that on the basis $a^l b^m$ of
$D_{q\pm}^2$,
where $m\in \N$, $l\in \Z$ and $a^{-l} := (a^*)^l$ for $l > 0$,
\begin{align}
\tau_1 (a^l b^m) &= \delta_{l}\delta_{m},
\nn\\
\tau_0^\up (a^l b^m) &=
\frac{1}{1\! -\! \lambda}\delta_{l}(1\! -\! \delta_{m}) -
\half\delta_{l}\delta_{m}
,
\label{eq:basis}\\
\tau_0^\dn (a^l b^m) &=
\frac{1}{1\! -\! \lambda}\delta_{l}(1\! -\! \delta_{m}) +
\half\delta_{l}\delta_{m}
\ ,
\nn
\end{align}
where $\lambda = (\pm q)^m$ and $\delta_k = 1$ when $k=0$ and $0$ otherwise.
Moreover (by checking on the basis)
\begin{align*}
\Tr_N(\pi_\pm (x))
&= (N + \sesq)\tau_1(x) + \tau_0^\up(x) + O(N^{-k})
\\
&= (N + \half)\tau_1(x) + \tau_0^\dn(x) + O(N^{-k})
\sepword{for all} x\in D_{q\pm}^2,  k > 0.
\end{align*}
For the following result, we shall use the following notation:\\
- denote the Wodzicki-type residue functional as in
\cite{ConnesMIndex}:
$$
\ncint T := \Res_{z=0} {\rm Tr}\, (T |D|^{-z}).
$$
- denote $r$ the projection onto the first two factors in
$\A(D^2_{q+}) \ox A(D^2_{q-}) \ox \A(S^1)$;
in particular, it yields a map
$$
r: \A(\Sf_q^*) \to \A(D^2_{q+}) \ox A(D^2_{q-}).
$$
- denote $T^0$ the grade-zero part
(with respect to the geodesic flow $\gamma (t)$ transported via
\eqref{eq:identific} to $\H'$)
of a diagonal operator $T$. 
~\\
A theorem in  \cite{DLSSV2}
states that the dimension spectrum of the spectral triple $(\A,\H,D)$ is
simple and given by $\{1,2,3\}$ with residues
\begin{align*}
\ncint T |D|^{-3} &= 2(\tau_1 \ox \tau_1) \bigl(r\rho(T)^0\bigr),
\\
\ncint T |D|^{-2} &= \bigl(\tau_1 \ox (\tau_0^\up + \tau_0^\dn)
+ (\tau_0^\up + \tau_0^\dn) \ox \tau_1\bigr) \bigl(r\rho(T)^0\bigr),
\\
\ncint T |D|^{-1}
&= (\tau_0^\up \ox \tau_0^\dn + \tau_0^\dn \ox \tau_0^\up)
\bigl(r\rho(T)^0\bigr).
\end{align*}
%
To see this recall that 
the dimension spectrum consists of the poles of the zeta function
$\zeta_T(z) := \Tr(T |D|^{-z})$ for all $T \in \Psi^0$,
but in our case, under the trace, it suffices to consider only 
$T\in P^\up\B$ or $T\in P^\dn \B$).\\
Moreover, in $\zeta_T(z)$ we can replace $T$
by $Q(\rho(T) \ox 1_2)Q$ since their difference is a smoothing
operator by \eqref{eq:symbol-smoothing}.
Calculating first for $P^\up T$ 
(fixing $s=\up$, splitting the overall trace into 
traces over $j, x, y$ and using the tracial property) yields,
up to holomorphic term,

$$
\Tr P^\up T |D|^{-z} =
\Tr(P^\up Q(\rho(T)  \ox  1_2)QP^\up\,|D|^{-z})
$$
$$
= \sum_{2j=0}^\infty (2j\!  +\!  \sesq)^{-z}
(\Tr_{2j} \ox \Tr_{2j\! +\! 1})
(r \rho(T)^0)
$$
$$
= (\tau_1 \ox \tau_1) (r \rho(T)^0)\,\zeta(z - 2) +
$$
\begin{equation}
 (\tau_1 \ox \tau_0^\dn + \tau_0^\up \ox \tau_1)
(r \rho(T)^0 )\,\zeta(z - 1) +
\label{eq:res-compute-up}
\end{equation}
$$
(\tau_0^\up \ox \tau_0^\dn) (r \rho(T)^0)\,\zeta(z). 
$$
Thus 
\begin{align}
\ncint P^\up T |D|^{-3} &= (\tau_1 \ox \tau_1) \bigl(r\rho(T)^0\bigr),
\nn \\
\ncint P^\up T |D|^{-2}
&= \bigl(\tau_1 \ox \tau_0^\dn + \tau_0^\up \ox \tau_1 \bigr)
\bigl(r\rho(T)^0\bigr),
\label{eq:expr-res-up} \\
\ncint P^\up T |D|^{-1}
&= (\tau_0^\up \ox \tau_0^\dn) \bigl(r\rho(T)^0\bigr).
\nn
\end{align}
Similar calculation for $P^\dn$, by shifting the summation index 
$j\mapsto j\! +\! \half$, 
yields just \eqref{eq:expr-res-up} with permuted tensor product.
Hence $\zeta_T$ has simple poles at $1$, $2$ and~$3$ with the residues
as stated.
%

\subsection{Local index formula}
Recall that with a general (odd) spectral triple $(\A,\H,D)$
the index of $D$ defines an additive map
\begin{equation}
K_1(\A) \to \Z,  ~~~
[U]\mapsto\ind(PUP),
\label{eq:fr-ind}
\end{equation}
where $U \in \Mat_r(\A)$ is a unitary representative of the $K_1$
class,
and $P=\half(1+F)$ with $F= D/|D|$
($PUP$ is automatically Fredholm).\\
This map
is computed by pairing $K_1(\A)$ with ``nonlocal'' cyclic
cocycles
$\chi_n$ given in terms of $F$
\begin{equation}
\chi_n(a_0, \dots, \a_n) = \lambda_n \Tr(a_0\,[F,a_1] \dots [F,a_n]),
\sepword{for all}  a_j \in \A,
\label{eq:nlcc}
\end{equation}
where $\lambda_n$ is a normalization constant and the
(smallest) integer $n \geq p$ is determined by the degree $p$ of
summability of the
Fredholm module $(\H,F)$ over $\A$.
In our case it is $1$-summable,
since the commutators $[F,\pi(x)]$, for $x\in\A$,
are trace-class (they are off-diagonal operators 
given by matrices of rapid decay).
Thus we need only the first Chern character
$\chi_1(a_0,a_1) =\lambda_1 \Tr(a_0\,[F,a_1])$, with $a_1,a_2 \in \A$.
An explicit expression for this cyclic cocycle on the basis 
of $SU_q(2)$ was obtained in~\cite{MasudaNW}.\\
~\\
On the other hand, the Connes--Moscovici local index theorem
expresses the index map in terms of a local cocycle $\phi_\odd$
in the $(b,B)$ bicomplex of $\A$, where
$$
b\varphi(a_0, a_1,\dots, a_{n+1}):=
$$
$$
\sum_{j=0}^n (-1)^j \varphi(a_0,\dots, a_j a_{j+1},\dots, a_{n+1})
+ (-1)^{n+1} \varphi(a_{n+1} a_0, a_1,\dots, a_n)
$$
and
$B = N B_0$,
with
\begin{align*}
& (N\psi)(a_0,\dots, a_{n-1}):= 
\sum_{j=0}^{n-1}  (-1)^{(n-1)j} \psi(a_j,\dots,a_{n-1}, a_0,\dots,a_{j-1})\ ,\\
 & 
(B_0 \varphi)(a_0,\dots, a_{n-1}):= 
\varphi(1,a_0,\dots, a_{n-1}) - (-1)^n \varphi(a_0,\dots,a_{n-1},1)\ ,
\end{align*}
satisfy $b^2 = 0$, $B^2 = 0$ and $b B + B b = 0$, so that $(b + B)^2 = 0$.\\
The cocycle $\phi_\odd$ is a local representative
of the cyclic cohomology class of $\chi_n$ (Chern character).
The pairing of the cyclic cohomology class
$[\phi_\odd] \in HC^\odd(\A)$ with $K_1(\A)$ gives
the index of $D$ with coefficients in $K_1(\A)$.
The (finite number of) components of $\phi_\odd = (\phi_1, \phi_3, \dots )$
are explicitly given in \cite{ConnesMIndex}
in terms of the operator $D$.\\
In our case
$(\A,\H,D)$ has metric dimension $3$ so 
the local cocycle $\phi_\odd$ has two components
$$
\phi_1(a_0,a_1) = \ncint a_0\, [D,a_1] \,|D|^{-1}
- \frac{1}{4} \ncint a_0\, \nabla([D,a_1]) \,|D|^{-3}
+ \frac{1}{8} \ncint a_0\, \nabla^2([D,a_1]) \,|D|^{-5},
$$
$$
\phi_3(a_0,a_1,a_2,a_3)
= \frac{1}{12} \ncint a_0\,[D,a_1]\,[D,a_2]\,[D,a_3]\,|D|^{-3},
$$
where $\nabla(T) := [D^2,T]$. 
(The cocycle condition 
$(b + B)\phi_\odd = 0$ reads $B \phi_1=0,~ b \phi_1 + B \phi_3 =0,~b\phi_3= 0$).\\
Whenever $[F,a]$ is traceclass  $\forall a \in \A$, 
which is our case, $\phi_1$, $\phi_3$ 
can be rewritten, using the binomial expansion, as
$$
\phi_1(a_0,a_1) = \ncint a_0 \,\delta(a_1) F|D|^{-1}
- \frac{1}{2} \ncint a_0 \,\delta^2(a_1) F|D|^{-2}
+ \frac{1}{4} \ncint a_0 \,\delta^3(a_1) F|D|^{-3},
$$
\begin{equation}
\phi_3(a_0,a_1,a_2,a_3) = \frac{1}{12} \ncint a_0 \,\delta(a_1)
\,\delta(a_2) \,\delta(a_3) F|D|^{-3}.
\label{eq:odd-cycle}
\end{equation}
~\\
Moreover Prop.~2 proved in \cite{ConnesSUq} applies to our case:
%
%
the Chern character $\chi_1$ is
equal to $\phi_\odd - (b + B) \phi_\ev$ where the cochain
$\phi_\ev = (\phi_0, \phi_2)$ is given by
\begin{align*}
\phi_0(a) &:= \Tr(Fa\,|D|^{-z}) \bigr|_{z=0},
\\
\phi_2(a_0,a_1,a_2)
&:= \frac{1}{24} \ncint a_0 \,\delta(a_1) \,\delta^2(a_2) F|D|^{-3}.
\end{align*}
\noindent
Remarks:\\\
- For the definition of $\phi_0$ it is necessary that $0\notin \Sigma$.\\
- The cochain $\phi_\ev = (\phi_0, \phi_2)$ is named $\eta$-cochain in
\cite{ConnesSUq}.\\
- In components, the equality of the characters means
$$
\phi_1 = \chi_1 + b \phi_0 + B \phi_2 , \qquad
\phi_3 = b \phi_2 .
$$
~\\
Furthermore there is a useful variant,
shown in \cite{DLSSV2} for $P$ of metric dimension 3, 
following Proposition 3 in \cite{ConnesSUq} 
valid for $P := \half(1 + F)$ of metric dimension 2:
the local Chern character
$\phi_\odd$ is equal to $\psi_1 - (b + B)\phi'_\ev$ where
$$
\psi_1(a_0,a_1) := 2 \ncint a_0 \,\delta(a_1) P|D|^{-1}
- \ncint a_0 \,\delta^2(a_1) P|D|^{-2}
+ \frac{2}{3} \ncint a_0 \,\delta^3(a_1) P|D|^{-3},
$$
and $\phi'_\ev = (\phi'_0, \phi'_2)$ is given by
\begin{align*}
\phi'_0(a) &:= \Tr(a\,|D|^{-z}) \bigr|_{z=0},
\\
\phi'_2(a_0, a_1, a_2)
&:= - \frac{1}{24} \ncint a_0 \,\delta(a_1) \,\delta^2(a_2) F|D|^{-3}.
\end{align*}
\noindent
Remarks:\\
- The term in $\psi_1$ involving $P |D|^{-3}$ is not present in \cite{ConnesSUq}
(as there the metric dimension of the projector $P$ is~$2$).\\
- Since $\phi_2 = - \phi'_2$ the last two propositions
jointly show that the cyclic $1$-cocycle $\chi_1$
can be given (up to coboundaries) in
terms of one single $(b,B)$-cocycle~$\psi_1$, more precisely
\begin{equation}
\chi_1 = \psi_1 -  b \beta,
\nn
\end{equation}
where
\begin{equation}
\label{eq:beta}
\beta (a) = 2\Tr(Pa\,|D|^{-z})\bigr|_{z=0}\ .
\end{equation}
~\\
%
%
%
To exemplify the theory above compute $\ind(PUP)$
when  $P=P^\up$ and $U$ is the unitary operator
\begin{equation}
U = \begin{pmatrix} a & b \\ -q b^* & a^* \end{pmatrix},
\label{eq:uni-def}
\end{equation}
acting on the doubled Hilbert space $\H \ox \C^2$ via the
representation $\pi \ox 1_2$.
Calculating $\psi_1(U^{-1}, U)$, with the local cyclic cocycle $\psi_1$ 
extended by $\Tr_\C^2$, yields
$$
 2  \sum_{kl} \ncint U_{kl}^* \,\delta(U_{lk}) P |D|^{-1}
- \sum_{kl}\ncint U_{kl}^* \,\delta^2(U_{lk}) P |D|^{-2}
+ \frac{2}{3} \sum_{kl} \ncint U_{kl}^* \,\delta^3(U_{lk}) P |D|^{-3}\ .
$$
Using \eqref{eq:expr-res-up}, simplifying due to 
$\rho(\delta^2(U_{kl})) = \rho(U_{kl})$, computing
(using the relations in $\A(D_{q\pm}^2)$)
$$
\rho(U_{kl}^*\,\delta(U_{lk}))^0 = 2(1 - q^2)\,1 \ox r_-(b)^2
$$
and substituting \eqref{eq:basis} produces
$$
\psi_1(U^{-1},U)
=
2(1 - q^2) (2\tau_0^\up \ox \tau_0^\dn + 2/3 \tau_1 \ox \tau_1)
\bigl(1 \ox r_-(b)^2 \bigr)
- (\tau_1 \ox \tau_0^\dn + \tau_0^\up \ox \tau_1)
\bigl(1 \ox 1 \bigr) = -2.
$$
Thus, including the correct normalization constant $-\half$,  
$$
\ind(PUP) = -\half \psi_1(U^{-1}, U) = 1.
$$

\section{Final comments}

Since $HC^1(\A(SU_q(2)))$ is one-dimensional, 
the characters of the 1-summable Fredholm modules found 
in~\cite{ConnesSUq} and in~\cite{DLSSV2}, 
are all cohomologous to the cyclic cocycle found in~\cite{MasudaNW}.\\
~\\
As signaled in the text, the polynomial algebras $\A$ employed above admit 
a smooth completion (consisting of series in the (ordered) generators 
with coefficients of rapid decay.
The fact that these pre C$^*$-algebras are e.g. closed under holomorphic 
functional calculus, follows from Lemma~2 of \cite{ConnesCIME} 
and the homomorphism properties of the maps \eqref{eq:q-disks}, 
\eqref{eq:symbol-map} and \eqref{eq:symbol-map-bis}.\\
~\\
It may be interesting to mention another class 
of infinitesimal operators with several interesting properties,
given by the two-sided ideal $\G$ of $\B(\H)$,
that was denoted by $\K_q$ in~\cite{Naiad}, \cite{DLSSV2}.
It is  generated by one operator $L_q$, where 
$L_q  v_{x,y,s}^j  := q^j \,v_{x,y,s}^j$,
but using the methods of \cite{DFWW} it can be seen that it is
in fact independent on the choice of $0<q<1$.
The ideal $\G$ is strictly contained in the ideal of infinitesimals
of arbitrary order, that is, compact operators whose $n$-th singular 
value $\mu_n$ satisfies $\mu_n = O(n^{-\a})$, for all $\a > 0$.
It seems however that it is too big for q-deformed groups,
for instance it is not so useful for the purpose
of \eqref{eq:symbol-map-bis} since the homomorphism
$\rho_\bullet$ does not extend to $\G$, even though
$\rho_\bullet (L_q) = 0$.\\
~\\
An interesting question is if the cochain $\beta$ \eqref{eq:beta} 
(whose coboundary is the difference 
between the original Chern character and the local one),
is given by the remainders in the rational approximation of the logarithmic
derivative of the Dedekind eta function as in ~\cite{ConnesSUq}.
The quite involved preliminary computations indicate that 
the higher derivatives are needed in the case at hand.\\
~\\
As far as the homogeneous spaces of $SU_q(2)$ are regarded, it should be added that 
local index formula on the equatorial Podle\'s quantum sphere $S_q^2$ 
has been worked out in \cite{lifeq} and in \cite{lifall} 
on the general family $S_{qc}^2$ of Podle\'s spheres.\\
~\\
A further study of these geometries and of $SU_q(2)$
presented in this note, which are characterized by $D$ with linearly
growing spectrum and the real structure and first order condition satisfied 
up to rapidly decaying operators, will be highly appreciated.
This in particular regards the investigation, along the lines of \cite{CoCh},
of the internal perturbations of the spectral action, which present
some difficulties (non-vanishing tadpole graph) for these examples.
\vspace{3pt}


\end{document}